\documentclass[11pt,a4paper]{amsart}

\usepackage{fancyhdr,stmaryrd}
\usepackage{amssymb} 
\usepackage{combelow} 
\usepackage{tensor} 
\usepackage{hyperref} 
\usepackage[all]{xy} 
\usepackage{mathrsfs} 
\usepackage{indentfirst} 
\usepackage{mathtools} 
\usepackage{titlesec} 
\usepackage[a4paper]{geometry}
\usepackage{xstring,xifthen}
\SetSymbolFont{stmry}{bold}{U}{stmry}{m}{n}
\usepackage{bm}
\usepackage{enumerate}
\usepackage{amsfonts,amsmath}

\newtheorem{theorem}{Theorem}[section]

\newtheorem{lemma}[theorem]{Lemma}

\theoremstyle{definition}

\newtheorem{remark}[theorem]{Remark}
\newtheorem{qu}[theorem]{Question}
\newtheorem{conj}[theorem]{Conjecture}

\newcommand{\HH}{\mathrm{HH}}
\newcommand{\Hc}{\mathrm{H}}
\newcommand{\hch}{\mathrm{hch.dim}}
\newcommand{\res}{\mathrm{res}}
\newcommand{\hh}{\mathrm{hh.dim}}
\newcommand{\gl}{\mathrm{gldim}}
\newcommand{\F}{\mathcal{F}}
\newcommand{\G}{\mathcal{G}}
\newcommand{\C}{\mathcal{C}}

\newcommand{\Aut}{\mathrm{Aut}}
\newcommand{\Hom}{\mathrm{Hom}}

\makeatletter
\newcommand{\@ndereq}[2]{%
  \vtop{
    \lineskiplimit\maxdimen
    \lineskip-.5\p@
    \ialign{$\m@th#1\hfil##\hfil$\crcr\sim\crcr#2\crcr}%
  }%
}
\newcommand{\sims}{\mathrel{\mathpalette\@ndereq{\scriptstyle \s}\relax}}
\makeatother

\newcommand{\titlename}	
{A reduction theorem  for non-vanishing of Hochschild cohomology of block algebras and Happel's property}
\newcommand{\shorttitlename}
{Reduction theorem non-vanishing  of Hochschild}

\newcommand{\authorname}      {Patrick Serwene$^1$ \and Constantin-Cosmin Todea$^2$}
\newcommand{\pdfauthorname}   {Patrick Serwene, Constantin-Cosmin Todea}
\newcommand{\shortauthorname} {P. Serwene \and C.-C. Todea}

\newcommand{\universitynameA}  {$^{1}$ Technische Universität Dresden}
\newcommand{\facultynameA}     {Faculty of Mathematics}
\newcommand{\departmentnameA}  {}
\newcommand{\addressA}  	   {01062 Dresden, Germany}

\newcommand{\universitynameB}  {$^{2}$Technical University of Cluj-Napoca}
\newcommand{\facultynameB}     {}
\newcommand{\departmentnameB}  {Department of Mathematics}
\newcommand{\addressB}  	   {G. Baritiu 25, RO-400027, Cluj-Napoca, Romania}

\newcommand{\emailaddressA}    {patrick.serwene@tu-dresden.de}
\newcommand{\emailaddressB}    {constantin.todea@math.utcluj.ro}

\newcommand{\articleabstract}{In this short research note we obtain a reduction theorem for the non-vanishing of the first Hochschild cohomology of  block algebras of finite groups with non-trivial defect groups. Along the way we investigate this problem for the blocks of some simple finite group algebras. Mimicking the case of  blocks of finite group algebras  we find some examples of category algebras that satisfy Happel's property.}

\newcommand{\msc}{20C20, }

\newcommand{\keywordterms}{finite, group, block, Hochschild, cohomology, Happel}

\def\depA{\departmentnameA}
\StrLen{\depA}[\depAlen]
\def\depB{\departmentnameB}
\StrLen{\depB}[\depBlen]
\def\facA{\facultynameA}
\StrLen{\facA}[\facAlen]
\def\facB{\facultynameB}
\StrLen{\facB}[\facBlen]

\newcommand{\institutionA}{
\universitynameA\\
\ifthenelse{\facAlen>0}{\facultynameA\\}{}
\ifthenelse{\depAlen>0}{\departmentnameA\\}{}
\addressA}
\newcommand{\institutionB}{
\universitynameB\\
\ifthenelse{\facBlen>0}{\facultynameB\\}{}
\ifthenelse{\depBlen>0}{\departmentnameB\\}{}
\addressB}
\hypersetup{colorlinks = true, linkcolor = blue, anchorcolor =red, citecolor = blue, filecolor = red, urlcolor = blue, pdfauthor=\pdfauthorname, pdftitle=\titlename, pdfkeywords=\keywordterms, pdfsubject=\articleabstract}
\titleformat{\section}{\Large\bfseries}{\thesection}{1em}{}
\titleformat{\subsection}{\large\it}{\thesubsection}{1em}{}
\fancypagestyle{mypagestyle}{
  \fancyhf{}
  \fancyhead[OC]{\textit{\shorttitlename}}
  \fancyhead[EC]{\textit{\shortauthorname}}
  \fancyfoot[C]{\medskip $\leftslice\,$\thepage$\,\rightslice$ }
  
}
\fancypagestyle{firstpagestyle}{
  \fancyhf{}
  \fancyfoot[C]{\medskip $\leftslice\,$\thepage$\,\rightslice$ }
  
}
\title[\shorttitlename]{\LARGE{\titlename}}
\author[\shortauthorname]{\large{\authorname}
\medskip\\
{\footnotesize \institutionA\medskip\\\institutionB\medskip\\
$\begin{array}{l}
\text{email}^1\text{: \texttt{\href{mailto:\emailaddressA}{\emailaddressA}}}\\
\text{email}^2\text{: \texttt{\href{mailto:\emailaddressB}{\emailaddressB}}}
\end{array}$
}}
\begin{document}
\begin{abstract}
\articleabstract\\[0.1cm]
\textsc{MSC 2010.} \msc\\[0.1cm]
\end{abstract}
\begingroup
\def\uppercasenonmath#1{} 
\let\MakeUppercase\relax 
\maketitle
\endgroup
\thispagestyle{firstpagestyle}


\section{Introduction} \label{sec1}
Let $k$ be a an algebraically closed field of characteristic $p$ and  $b$ a block idempotent of $kG,$ where $G$ is a finite group such that $p$ divides the order of $G.$  It is known that the first Hochschild cohomology group (denoted by $\HH^1(kG)$) of the group algebra $kG$  forms a Lie algebra over $k$ and there is recent growing interest on the solvability of first Hochschild cohomology group of finite dimensional algebras in general \cite{RuSol} and block algebras  of finite groups in particular \cite{LiRu}. We shall use the following notation
\[\G:=\{(G,P) \mid G \text{ is a finite group s.t.}\  p \mid |G|, P\text{ a non-trivial } p\text{-subgroup of } G\}.\]
The non-vanishing  of $\HH^1(kG)$ is obtained as a consequence of  a result of Fleischmann, Janiszczak and Lempken, \cite{FJL}. We know that the block algebra $kGb$ is a simple $k$-algebra if and only if $b$ has trivial defect group. However, it is still an open problem if  $\HH^1(kGb)\neq 0$ for any block algebra $kGb$ with non-trivial defect group. The following question was launched by Linckelmann at ICRA 2016, Syracuse.
\begin{qu}\label{questLi}(\cite[Question 7.7]{Li}) Is it true that for any $G$ such that $(G,P)\in\G$ and any block algebra  $kGb$ with defect group $P$
we have $\HH^1(kGb)\neq 0$?
\end{qu}

We will investigate the so-called Happel's property for Hochschild cohomology and we spend some time to recall the concepts. For a finite dimensional $k$-algebra $A$, we denote the global dimension of $A$ as usual by $\gl(A)$ and the Hochschild cohomological (respectively, homological) dimension of $A$ by
$$\hch(A):=\mathrm{sup}\{n\in\mathbb{N} \mid \HH^n(A)\neq 0\},$$
(respectively, $\hh(A):=\mathrm{sup}\{n\in\mathbb{N} \mid \HH_n(A)\neq 0\}$). We follow and are inspired by the survey \cite{Cr}.
The statement that for any finite dimensional $k$-algebra, $\hch(A)$ is finite implies that $\gl(A)$ is finite, is called \textit{Happel's property}. The  analogous statement for Hochschild homology (i.e. $\hh(A)<\infty$ implies $\gl(A)<\infty$) is called Han's property. Inspired by block algebras of finite groups and group algebra cases,  we are interested in studying Happel's property for some Frobenius $k$-algebras. The converse statement to Happel's property is obtained by theorem due to Keller \cite{Kel} (see \cite[Theorem 2.1]{HapZac}), given in the case of Hochschild homology for any finite dimensional $k$-algebra.

In the following paragraph we include remarks and questions suggested by an email discussion with Markus Linckelmann.  We assume that $b$ has  non-trivial defect group $P$, hence the block algebra $kGb$ is not semisimple and $\gl(kGb)=\infty.$  So, Happel's property for $kGb$ is the same to showing  that $\hch(kGb)=\infty,$ (equivalently, $\hch(kGb)\neq 0$.)  The Hochschild cohomology $\HH^*(kGb)$ is always infinite-dimensional, because its Krull dimension is positive, equal to the rank of $P$,
so $\HH^*(kGb)$ contains a polynomial subalgebra. One can say a bit more (see \cite[Theorem 13.1]{Licra10}):
the longest sequence of consecutive integers $n$ for which $\HH^n(kGb)=0$ is bounded in terms
of $P.$ In particular, Happel's property is true for $kGb.$ A valid question to ask would be this: what can one say about the smallest positive integer n such
that $\HH^n(kGb)$ is non-zero? By the above, we know that $n$ is bounded in terms of $P.$ And as we can see in Question \ref{questLi}, we expect that $n=1.$

As an alternative method of proof for block algebras, Happel's property is a consequence of the validity of a well-known conjecture in the theory of saturated fusion systems, Conjecture \ref{conj13}; for this theory we follow \cite{AKO}. A saturated fusion system $\F$ on a finite $p$-group $S$ is called  \textit{realizable} if there is a finite group $G_1$ such that $\F=\F_{S}(G_1),$ where $S$ is a Sylow $p$-subgroup of $G_1;$ otherwise $\F$ is said to be exotic.  We consider a $p$-subgroup $P$ of $G$ as a defect group of $b.$ We fix a maximal $b$-Brauer $(P,e)$ and denote by $\F$ the saturated fusion system $\F_{(P,e)}(G,b)$ of the block $b.$ This means that $b$ is an $\F$-block and we say for shortness that $kGb$ is an $\F$-\textit{reduction simple} block algebra if $P$ has no non-trivial proper strongly $\F$-closed subgroups. We restate the following equivalent form of  \cite[Part 4, Section 7, p.298]{AKO}:
\begin{conj}\label{conj13}(\cite[Conjecture 1.1]{serwene}) \textit{Let $\F$ be a saturated fusion system. For any finite group $G$ having an $\F$-block $b$, it follows that $\F$ is realizable}.  
\end{conj}
In the next remark  we present the ingredients that assure us that if Conjecture \ref{conj13} is verified for $b$ then any $kGb$ satisfies Happel's property.
\begin{remark}\label{rem13} Let $(G,P)\in \G,$  a saturated fusion system $\F$ on $P$ and $b$   an $\F$-block. If $\F$ is realizable it follows that there is a finite group $G_1$ with $P$ a Sylow $p$-subgroup such that $\F=\F_P(G_1).$ Using the Cartan-Eilenberg Stable Elements Theorem \cite[Theorem 10.1]{CarEil}, we know that the cohomology of $\F$ is
$(\Hc^n(P,k))^{\F}\cong \Hc^n(G_1,k),$
for any positive integer $n.$
Recall that the cohomology of the block $b$ is $(\Hc^n(P,k))^{\F},$ that was introduced by Linckemann in \cite{LiTr} and that there is an injective map
\begin{equation}\label{eq2}
(\Hc^n(P,k))^{\F}\rightarrow \HH^n(kGb),
\end{equation}
for any positive integer $n, $ see \cite[Theorem 5.6 (iii)]{LiTr}. By \cite[Corollary 1]{Swan}, since $p$ divides the order of $G$, it follows that $\Hc^n(G_1,k)$ is non-zero for infinitely many values of $n>0,$ hence applying (\ref{eq2}) we obtain that $\hch(kGb)=\infty.$

\end{remark}
In the first main theorem we present some Frobenius $k$-algebras (similar to group algebras) that verify Happel's property. Let $\C$ be a finite category  (the class of morphisms is a finite set) and $k\C$ the category algebra. Let $\underline{k}$ be the constant functor  from $\C$ to the category of vector $k$-spaces, sending every object to $k$ and every morphism to identity of $k.$ The ordinary cohomology of $\C$ with coefficients in $k,$ denoted by $\Hc^*(\C, \underline{k}),$ can be defined as $\mathrm{Ext}^*_{k\C}(\underline{k},\underline{k}),$ see \cite{Xu}. We denote by $\bold{G},$ the group $G$ viewed as a category with one object, homomorphisms the elements of $G$ and the composition is the multiplication in $G.$ If $\pi:\C\rightarrow \bold{G} $ is a covariant functor, there is always a restriction map $\res_{\pi}^*:\Hc^*(\bold{G},\underline{k})\rightarrow \Hc^*(\C,\underline{k}),$ see \cite{Xu, Xu2}.
\begin{theorem}\label{thmcross}
Let $\C$ be a finite category such that $k\C$ is a Frobenius $k$-algebra. If there is a coovariant functor $\pi: \C\rightarrow \bold{G}$ such that $\res_{\pi}^*:\Hc^*(\bold{G},\underline{k})\rightarrow \Hc^*(\C,\underline{k})$ is injective, then $k\C$ satisfy Happel's property. 
\end{theorem}
We do not claim that in the above theorem we obtain remarkable new examples, but we just want to emphasize the method described in Remark \ref{rem13}. It may happen that the validity of Happel's property for the examples in Remark \ref{rem21} (that verify the assumptions of Theorem \ref{thmcross}) can be obtained by using other more simple arguments.

In the second theorem of this paper we return to Question \ref{questLi} and we obtain the following reduction result: if for any $(L,Q)\in \G$  with $L$ quasi-simple, for any block algebra $kLd$ with defect group $Q,$ we have $\HH^1(kLd)\neq 0$  then for any $(G,P)\in \G,$ with $\Aut(P)$ a $p$-group and for any $\F$-reduction simple block algebra $kGb$ with defect group $P,$ we have $\HH^1(kGb)\neq 0.$ 
 
\begin{theorem}\label{thm:13}
If there is $(G,P)\in \G,$ with $\Aut(P)$ a $p$-group and there is an $\F$-reduction simple block algebra $kGb$ with defect group $P,$  such that  $\HH^1(kGb)=0$  then there is $(L,Q)\in \G$  with $L$ quasi-simple and there is a block algebra $kLd$ with defect group $Q$ such that $\HH^1(kLd)= 0.$
\end{theorem}

In the final theorem we add, to the list of Murphy \cite{Mur}, of  Benson, Kessar and Linckelmann \cite{bkl} and Todea \cite{Tod}, new classes of block algebras for which the answer to Question \ref{questLi} is true.
For the first statement of the next theorem let $\tilde G=S_n, G=A_n$ and we assume that $p$ is odd. Recall that $G \unlhd \tilde G$ with $|\tilde G/G|=2$, hence $\tilde G/G \cong C_2$, which is a $p'$-group. We consider the cyclic group $C_2=<\sigma>,$ generated by $\sigma.$ Let $\tilde b$ be a block of $k \tilde G$ covering $b$. Clearly $\tilde b$ has the same defect group as $b$, see for example \cite[Proposition 2.22]{serwene2} and it is known that the block $b$ is $C_2$-stable. Let $(K,\mathcal{O}, k)$ be a $p$-modular system such that $K$ contains $\varepsilon,$ a primitive second root of unity. The group of characters of $C_2,$ that is $\widehat{C_2}:=\Hom(C_2,K^\times),$ is isomorphic with $C_2$ and we have that $\widehat{C_2}=<\widehat{\sigma}>,$ where $\widehat{\sigma}(\sigma)=\varepsilon.$  This group acts on the blocks of $k \tilde G$ and we adopt the notation $^{\widehat{\sigma}} \tilde b=b^\ast,$  hence $b^\ast$ is a new  block of $k \tilde G$.

\begin{theorem}\label{thm16} Let $G$ be a finite group, $p$ be  a prime dividing $|G|$ and $b$ a block of $kG$ with non-trivial defect group. We assume that one of the following statements  is true.

\begin{itemize}

\item[(i)]  $\tilde G=S_n, G=A_n,$ $p$ is odd  and $\tilde b\neq b^*;$ 
\item[(ii)] $G=J_1,$ the first Janko group;
\item[(iii)] $G=J_2,$ the Hall-Janko group;
\item[(iv)] $G=J_4,$ the fourth Janko group, but $b$ is not the principal $3$-block and it is not the unique $3$-block of maximal defect group $3_{+}^{1+2}$.
\end{itemize}
Then $\HH^1(kGb) \neq 0.$ 
\end{theorem}
Section \ref{sec2} contains the proofs of the above theorems.

\section{Proof of the Theorems} \label{sec2}
\begin{proof}\textbf{(of Theorem \ref{thmcross})}
Since $k\C$ is a Frobenius $k$-algebra it follows by \cite[Exercise 4.2.2]{Wei}, that $\gl(k\C)$ is zero or infinity.
But  $k\C$ cannot be semisimple, because  in this case $\hch(k\C)$ should be zero and we will show  that $\hch(k\C)$ is infinity, using our assumption.

We know that $\Hc^*(\bold{G},\underline{k})\cong \Hc^*(G,k)$ and since $p$ divides $|G|$ it follows by \cite[Corollary 1]{Swan} that $\Hc^n(\bold{G},\underline{k})$ is non-zero for infinitely many values of $n>0.$ Thus, $\Hc^n(\C,\underline{k})$ is is non-zero for infinitely many values of $n>0,$ since $\res_{\pi}^*$ is injective.

We conclude by applying \cite[Theorem A]{Xu} to get that   $\Hc^n(\C,\underline{k})$ is a direct summand of $\HH^n(k\C)$ (for any non-negative integer $n$) hence $\hch(k\C)=\infty.$

\end{proof}
\begin{remark}\label{rem21}

1) Let $\mathcal{P}$ be a finite $G$-poset such that $\mathcal{P}$ has only the equality relation between elements (any two different elements are incomparable). The category algebra $k(G\propto \mathcal{P})$ (over the transporter category $G\propto \mathcal{P}$, see \cite[2.1]{Xu2}) is a Frobenius $k$-algebra. It is in fact a symmetric $k$-algebra, since it is a groupoid algebra and we can apply \cite[Theorem 1.1]{LiInv}. If the Euler characteristic
$\chi(\mathcal{P}, k)$ (\cite[2.2]{Xu2}) is invertible in $k$ (in this case $\chi(\mathcal{P}, k)=|\mathcal{P}|$) then, there is an injective restriction map
$\Hc^*(\bold{G}, \underline{k})\rightarrow \Hc^*(G\propto \mathcal{P},\underline{k}),$
see \cite[Theorem 1.1 ]{Xu2}.

2) A category algebra $k\C$ over a finite inverse category $\C$ is a symmetric $k$-algebra \cite[Theorem 1.1]{LiInv}. We can look for more examples of finite inverse category algebras satisfying the assumptions of Theorem \ref{thm:13} by using \cite[Theorem 1.5]{LiInv}.
\end{remark}

The following lemma is folklore.
\begin{lemma}\label{lemma:22}
Let $G, H$ be any finite groups (with $p$ dividing their orders) and $kG=\bigoplus_i kGb_i,$ $kH=\bigoplus_j kHb'_j$ the block decompositions of their respective group algebras. Then the block decomposition of $k(G \times H)$ is $k(G \times H) \cong \bigoplus_{i,j} kGb_i \otimes kHb'_j.$
\end{lemma}
The proof of Theorem \ref{thm:13} is an almost verbatim translation of proof given by Kessar for \cite[Theorem 3.1]{art:Kess2006}. Since some parts  are relatively different we explicitly present the proof for the convenience of the reader.
\begin{proof}\textbf{(of Theorem \ref{thm:13}).}

We assume that $(G,P)\in\G$ with $\Aut(P)$ a $p$-group and with $G$ of minimal order with respect to having an $\F$-reduction simple block $b$ of defect $P,$  such that $\HH^1(kGb)=0.$ We split the  proof in several steps.

\textit{Step 1.} For any $H,$ a normal subgroup of $G,$ such that $p$ divides the order of $H$ and such that $c$ is a block of $kH$ covered by $b,$ it follows that $c$ is $G$-stable. \\
This is true by a well-known Clifford theoretic fact and by the minimality of $|G|,$ see \cite[Proposition 2.13]{art:Kess2006}.

\textit{Step 2.} The group $G$ is generated by all conjugates of $P.$ \\ For showing this, let $H:=<{}^gP| g\in G>,$ which  is a normal subgroup of $G.$ Using Step $1$, let $c$ be the unique block of $kH$ covered by $b.$ We know that $G$ acts by conjugation on $kHc$ and hence we  have a group homomorphism from $G$ to the group of exterior automorphisms (as $k$-algebras) of $kHc.$  We denote by $K$ the normal subgroup of $G$ obtined as the kernel of this group homomorphism, see \cite[3.2]{art:Kess2006} or \cite[Section 5]{art:Kuels1995}. By Step 1 we consider $c_1$ the unique block of $kK$ covered by $b.$ By \cite[Theorem, p.303 ]{art:Kuels1995} we know that $G/K$ is a $p'$-group whose order is dividing $|\mathrm{Out}(P)|^2.$ Moreover $c_1$ is in fact $b$ (having the same defect group $P$ in $K$) which remains a block of $kK$ and  the block algebra $kGb$ is a  crossed product  of $G/K$ with $kKb.$ The minimality assumption give us $G=K.$ Since $G=K$ and $b,c$  are blocks of $kG,kH,$ respectively, having a common defect group $P,$ by \cite[Theorem 7]{art:Kuels1990} we get that $kGb$ and $kHc$ have isomorphic source algebras. It follows that $$\HH^1(kHc)\cong \HH^1(kGb)=0$$ hence again the minimality assumption forces the equality $G=H.$

\textit{Step 3.} For any proper normal subgroup $N$ of $G,$ if $d$ is a block of $kN$ covered by $b,$ then $d$ is of defect zero.\\
By Step 1 we assume that $d$ is $G$-stable and then by \cite[Theorem 6.8.9, (ii)]{LiBookII} we have that $P\cap N$ is a defect group of $d.$ If $P\cap N\neq \{1\},$ since $P\cap N$ is always strongly $\F$-closed and $b$ is $\F$-reduction simple, it follows that $P\cap N=P.$ That is $P\subseteq N,$ hence $P$ and all its conjugates are in $N.$ Using Step 2 it follows that $N=G,$ a contradiction.

\textit{Step 4.} Under  the assumptions of Step 3, by a variation of second Fong's reduction (see \cite[Proof of Theorem 3.1]{art:Kess2006}), there is a central $p'$-extension
\[
\xymatrix{
1\ar[r] & Z\ar[r]&\tilde{G}\ar[r]& \bar{G}\ar[r]&1,}
\]
where $\bar{G}$ is $G/N$ and there is a block $c$ of $k\tilde{G}$ such that $kGb$ is Morita equivalent with $k\tilde{G}c.$ In particular we obtain  $$\HH^1(k\tilde{G}c)\cong \HH^1(kGb)=0.$$

\textit{Step 5.} We assume now that $N$ is maximal proper subgroup of $G.$ Applying Step 4 in this case, we get a central $p'$-extension

\[
\xymatrix{
1\ar[r] & Z\ar[r]&\tilde{G}\ar[r]& \bar{G}\ar[r]&1,}
\]
where $\bar{G}:=G/N$ is a simple finite group. Note that $Z(\tilde G)/Z \unlhd \tilde G/Z \cong \bar{G},$ hence $Z=Z(\tilde G)$ is a cyclic $p'$-subgroup. 
Again by Step 4, there is an $\mathcal F$-block $c$ of $k \tilde G$ such that $\HH^1(k \tilde G c)=0$ and $\HH^2(k \tilde G c)=0$. Let $L=[\tilde G, \tilde G]$.\\
First note that if $\tilde G/Z \cong \bar{G}$ is a cyclic group, then $\tilde G$ must be abelian, see proof of \cite[Lemma 3.4]{serwene}, hence $L=1$ in that case. Since $ \tilde{G}$ abelian it follows that $Z=\tilde{G},$ hence $\bar{G}=1,$ a contradiction.\\
We are left with the case that $\bar{G}$ is non-abelian simple. We have that $\tilde G=LZ$, $L$ is a quasi-simple group and $Z(L) \subseteq Z $ is a cyclic $p'$-group, see again proof of \cite[Lemma 3.4]{serwene}. We obtain another short exact sequence:

\[
\xymatrix{
1\ar[r] & K\ar[r]& L \times Z \ar[r]^{\pi} & \tilde{G} \ar[r]&1,}
\]
where $\pi (e,z)=ez$, for any $e\in L,z\in Z$ and $$K:=\ker \pi=\{ (e,e^{-1}) \mid e \in L \cap Z \}.$$
By \cite[Lemma 2.10(i)]{fk}, for a block $c$ of $k\tilde{G}$, there is a unique block $c'$ of $k(L \times Z)$ such that $c$ is dominated by $c'$, that is $\pi(c') \neq 0;$ where $\pi(c')$ is an idempotent of $Z(k \tilde G)$ and $\pi$ is the $k$-algebra surjection $\pi: k(L \times Z) \rightarrow k \tilde G$ (by abuse of notation) induced by the canonical surjection $\pi:L \times Z \rightarrow \tilde G,$ where $\tilde G \cong (L \times Z)/K.$ Since $K$ is a $p'$-subgroup and $\pi(c') \neq 0,$ by \cite[Lemma 2.10(ii)]{fk} it follows that $\pi(c')=c$ and $k \tilde G c \cong k(L \times Z)c'$ as $k$-algebras. Thus, we obtained a block $c'$ of $k(L \times Z)$ such that $\HH^1(k(L \times Z)c')=0.$ 

For the last part, utilizing Lemma \ref{lemma:22}, it is well known that if $c'$ is a block of $k(L \times Z),$ then there is a block $d$ of $kL$ and a block $d'$ of $kZ$ (hence of defect zero) such that $$k(L \times Z)c' \cong kLd \otimes kZd'.$$
It follows that: $$0=\HH^1(k(L\times Z)c' \cong \HH^1(kLd \otimes kZd') \cong (\HH^1(kLd) \otimes Z(kZd')) \oplus (Z(kLd) \otimes \HH^1(kZd')$$
$$=(\HH^1(kLd) \otimes k) \oplus (Z(kLd) \otimes 0)=\HH^1(kLd),$$ hence $d$ is a block of $kL$ with $\HH^1(kLd)=0.$ 
\end{proof}
It remains a challenge (for future work) to remove the assumptions that $\Aut(P)$ is a $p$-group and that the block $b$ is $\F$-reduction simple.
\begin{proof}\textbf{(of Theorem \ref{thm16}).}  

(i).  If $\tilde b \neq b^\ast,$ then the block algebras $kS_n \tilde b$ and $kA_nb$ are isomorphic, hence $$\HH^1(kA_nb)\cong \HH^1(kS_n \tilde b) \neq 0,$$ where the last statement is true by \cite[Theorem 1.1]{bkl}.

(ii). $|J_1|=175560=2^3 \cdot 3 \cdot 5 \cdot 7 \cdot 11 \cdot 19$.\\

\textit{Case (1)}: $p$ odd - then all blocks $b$ that have a non-trivial defect group must have a cyclic defect group. In this case we know that $\HH^1(kGb) \neq 0$.\\

\textit{Case (2)}: $p=2$ - using the character table from \cite{Janko} and the introduction from \cite{LandMich}, it is easy to see that all $2$-blocks of $J_1$ have defect zero or one, except for the principal $2$-block whose defect group is $C_2 \times C_2 \times C_2$. It follows that a $2$-block $b$ of $J_1$ with non-trivial defect group has either a cyclic defect group, in which case we know that $\HH^1(kGb) \neq 0$ or its defect group is $C_2 \times C_2 \times C_2$ and $b$ is principal. For the second case we apply Theorem 6.10 of \cite{Rou} to obtain that Brou\'e's Abelian Defect conjecture is true and get $\HH^1(kGb) \neq 0$.\\

(iii). $|J_2|=604800=2^7 \cdot 3^3 \cdot 5^2 \cdot 7$. We shall use the description of the blocks of this group that appears  in \cite{sin}.\\

If $p=7$, we have a cyclic defect group and we are done.\\

Let $p=5$, then there are five $5$-blocks. Three of these have defect zero and one has defect one. For those four blocks there is nothing to prove. The principal $5$-block has defect group $C_5 \times C_5$. By \cite{Holloway} we see that Brou\'e's Abelian Defect Group Conjecture holds in this case and we are also done with this block.\\

Now assume $p=3$. Then there are three blocks each of defect zero and  three blocks of defect $1$, for which there is nothing to check. The only other block is the principal $3$-block, which we denote by $b_0$ and $kJ_1b_0=B_0$. We shall use the strategy developed in \cite{Mur} to obtain $\dim_k \HH^1(B_0)$. For this, we use GAP \cite{gap} to obtain that $\dim_k \HH^1(kG)=7$. We have a decomposition $$kG=B_0 \times B_1 \times  B_2 \times B_3 \times B_4' \times B_5' \times B_6',$$ where $B_i$ for $i \in \{1,2,3\}$ are the block algebras of defect one and $B_i'$ for $i \in \{4,5,6\}$ are the block algebras of defect zero. Then $$\HH^1(kG)=\HH^1(B_0) \oplus \HH^1(B_1) \oplus \HH^1(B_2) \oplus \HH^1(B_3),$$ hence 
\begin{equation}\label{eq7}
\dim_k \HH^1(B_0)=7-\sum\limits_{i=1}^3 \dim_k \HH^1(B_i).
\end{equation}
Let $i \in \{1,2,3\}$. It is well known that $\dim_k \HH^1(B_i)=\frac{|P|-1}{|E_i|}$, where $E_i$ is the inertial quotient with order dividing $2$. One can use GAP \cite{gap} to be more precise and obtain $|E_i|=2$, hence $\dim_k \HH^1(B_i) = 1$. It follows that 
\begin{equation}\label{eq8}
\sum\limits_{i=1}^3 \dim_k \HH^1(B_i) =3.
\end{equation}
By (\ref{eq7}), (\ref{eq8}) we obtain that $\dim_k \HH^1(B_0)=4$ as desired.\\

Finally, assume $p=2$. There are only two blocks, one has defect group isomorphic to $C_2 \times C_2$, for which Brou\'e's Abelian Defect Group Conjecture is known and the principal $2$-block. Using GAP \cite{gap}, we note that $\dim_k \HH^1(kG)=17$. As above, if $$kG=B_0 \times B_1, B_0=kGb_0, B_1=kGb_1,$$ then
\begin{equation}\label{eq9}
\dim_k \HH^1(B_0)=17-\dim_k \HH^1(B_1). 
\end{equation}
Since $B_1$ is a block algebra with defect group the Klein four group $l(B_1) \in \{1,3\}$. In any case, by \cite[Theorem 2.1]{holm}, we get that $\dim_k \HH^1(B_1) \in \{8,2\}$, thus
\begin{equation}\label{eq10}
\dim_k \HH^1(B_1) \leq 8
\end{equation}
and finally using (\ref{eq9}) and (\ref{eq10}), we get that $\dim_k \HH^1(B_0)>0$.

(iv). $|J_4|=2^{21}\cdot 3^3\cdot 5\cdot 7\cdot 11^3\cdot 23\cdot 29\cdot 31\cdot 37\cdot 43.$ We shall use the descriptions of the blocks (with non-trivial defect group) of this group that appears  in \cite[Lemma 5.1]{An} and \cite[Lemma 3.2]{Kos}. We may assume that $p\in\{2,3,11\},$ since any block of $kG$ for the other primes has cyclic or trivial defect group.

For $p=2$ the principal block $b_0$ of $kG$ is the unique block of positive defect (see \cite[Lemma 5.1 (b)]{An}), hence $$\HH^1(kGb_0)\cong \HH^1(kG)\neq 0.$$

For $p=3$ \cite[Lemma 5.1 (a)]{An} assure us that there are seven blocks of positive defect, giving the decomposition
$$\HH^1(kG)=\HH^1(B_0) \oplus \HH^1(B_1) \oplus \HH^1(B_2) \oplus \HH^1(B_3)\oplus \HH^1(B_4)\oplus \HH^1(B_5)\oplus \HH^1(B_6),$$
The block algebra $B_2$ is the unique $3$-block algebra of defect group an elementary abelian $3$-group of order $9.$ By \cite[Theorem 1.3 ]{Kos} Brou\'{e}'s Abelian defect group conjecture is true in this case, hence $\HH^1(B_2)\neq 0.$ Moreover, $\dim_k(\HH^1(B_2))=\dim_k (\HH^1(kA_8b_0'))$ (see \cite[Theorem 1.4]{Kos}), where $b_0'$ is the principal block of $kA_8.$ All the other block algebras $B_i,i\in\{3,4,5,6\}$ have cyclic defect group $C_3.$

For $p=11$ we know by \cite[Lemma 3.2]{Kos} that the principal block is the unique block of positive defect, hence  has its first Hochschild cohomology isomorphic with the non-vaninshing first Hochschild cohomology of the entire group algebra $kG;$ see \cite[Table 6, p.464]{Bla}.
\end{proof}

In the above proof (case $G=J_4,p=3$) it remains a challenge to deal with the non-vanishing of the first Hochschild cohomomology of the principal block algebra $B_0$ and of $B_1,$ the block algebra of maximal defect group $ 3_+^{1+2}.$

\textbf{Acknowledgments.} We thank Professor Shigeo Koshitani for suggesting the  reference \cite{sin}.

\phantomsection

\end{document}